\documentclass[reqno, 12pt]{amsart}

\usepackage{color}

\newtheorem{theorem}{Theorem}[section]
\newtheorem{corollary}[theorem]{Corollary}
\newtheorem{lemma}[theorem]{Lemma}

\theoremstyle{definition}
\newtheorem{remark}[theorem]{Remark}

\theoremstyle{definition}

\theoremstyle{definition}
\newtheorem{assumption}[theorem]{Assumption}

\numberwithin{equation}{section}

\makeatletter
\def\dashint{\operatorname%
{\,\,\text{\bf--}\kern-.98em\DOTSI\intop\ilimits@\!\!}}
\makeatother

\newcommand{\WO}[2]{\overset{\scriptscriptstyle0}{W}\,\!^{#1}_{#2}}

\def\coa{{\sf a}}
\def\cob{{\sf b}}
\def\coc{{\sf c}}

\def\bR{\mathbb{R}}

\def\bQ{\mathbb{Q}}

\def\cA{\mathcal{A}}
\def\cB{\mathcal{B}}
\def\cC{\mathcal{C}}

\title {Parabolic equations with measurable coefficients}
\author{Doyoon Kim \and N. V. Krylov}

\address{127 Vincent Hall, University of Minnesota, Minneapolis, MN
55455} 
\email{dykim@math.umn.edu}

\address{127 Vincent Hall, University of Minnesota, Minneapolis, MN
55455} 
\email{krylov@math.umn.edu}

\thanks{The second author was partially supported by
NSF Grant DMS-0140405}



\begin{document}

\begin{abstract} We investigate the unique solvability of second order
parabolic equations in non-divergence form in $W_p^{1,2}((0,T) \times
\bR^d)$, $p \ge 2$.  The leading coefficients are only measurable in
either one spatial variable or time and one spatial variable.  In
addition, they are VMO (vanishing mean oscillation) with respect to the
remaining variables.  
\end{abstract}

\maketitle

\section{Introduction}

This paper is a natural continuation 
of our  previous investigations \cite{Kr05}, \cite{KiK}.
By combining the techniques from these articles we investigate
parabolic equations of type
\begin{equation}                              \label{para_main_eq} 
u_t + a^{jk}(t,x) u_{x^j x^k} +
b^{j}(t,x) u_{x^j} + c(t,x) u = f
\end{equation} 
in Sobolev spaces $W^{1,2}_{p}$ 
with $p\geq2$ and the coefficients   being
just measurable in $x^{1}$ but VMO with respect to other variables.
Here 
$$
(t,x)\in\bR^{d+1}=\{(t,x^{1},x'):
t,x^{1}\in\bR,x'\in\bR^{d-1}\}
$$
and the equation is assumed to be uniformly nondegenerate with bounded
coefficients.

One of the advantages of having a ``good" theory for
such equations is demonstrated in \cite{KiK}
while treating the Dirichlet and Neumann problems,
the issues addressed in this paper as well.

The  amazing fact that there is a solvability
theory in Sobolev spaces for
 elliptic and parabolic equations
with discontinuous but VMO coefficients was discovered 
in \cite{CFL1}, \cite{CFL2}, and \cite{BC}. Before that
the Sobolev space theory was established for
some other types of discontinuities \cite{Lo}, \cite{Lo1},
\cite{Ch}, \cite{Sa} (see also \cite{Ki}  (\cite{Ki1}) for a modern approach covering $p\ne2$
 in the elliptic  (parabolic) case).
Solvability theory for discontinuous coefficients
is important not only from pure theoretical point of
view but also from the point of view of applications,
for instance, to random diffusion processes, see, for instance,
\cite{SV}, \cite{Kr05}. Observe that the class
of equations with VMO coefficients and
the class
of equations with discontinuities
 treated in \cite{Lo},  
\cite{Ch}, \cite{Sa}, \cite{Ki}, \cite{Ki1} have no common members
apart from the equations with just continuous coefficients.
In this paper we show that there is a unified approach
to both cases allowing one to treat equations
possessing one properties with respect to some variables
and other properties with respect to  the remaining ones.
Here we show that the coefficients $a^{jk}(t,x^{1},x')$
 may be just measurable in $x^{1}$ and VMO in $(t,x')$.
Even though the equations (and partly the results) of the present article
and \cite{KiK} are more general than those from
\cite{BC}, \cite{Ch}, \cite{Ki1}, \cite{Sa}, they are not 
general enough  to absorb \cite{Kr05}, where
equations are considered whose coefficients $a^{jk}$ are allowed to be 
{\em measurable\/} in $t$
and VMO in $x$. 
Furthermore,  the results here cover those of \cite{Ki1}
only for $p\geq2$. 
On the other hand, in \cite{Kr05}
and \cite{Ki1}
the coefficients only measurable in $x^{1}$ are not allowed.
Thus, the classes of equations here and in \cite{Kr05},
\cite{Ki1} are quite different.

It is worth noting that after \cite{CFL1}, \cite{CFL2},
\cite{BC} there were very many publications on
elliptic and parabolic equations with VMO coefficients
(see, for instance, the above mentioned references and
 \cite{HHH}, \cite{MPS}, \cite{PS},
\cite{So1}, \cite{So2}, and many references therein).
The approach we employ here is quite different from
the approaches of other authors and is taken from~\cite{Kr05}.

This paper is organized as follows. In section \ref{main_sec} we present
our main results. The case $p = 2$ is investigated in section
\ref{L_2_case}. In section \ref{aux_W_p} we present some auxiliary
results which are needed for the proof of Theorem \ref{para_main_01}. In
section \ref{proof_th_p} we prove Theorem \ref{para_main_01}.

A few words about notation.  As is seen from the above
 by $(t,x)$ we denote a point in $\bR^{d+1}$,
i.e.,
$(t,x) = (t,x^1,x') \in \bR \times \bR^{d} = \bR^{d+1}$, where $t \in
\bR$, $x^1 \in \bR$, $x' \in \bR^{d-1}$, and $x = (x^1,x') \in \bR^d$. By
$|u|_{0}$ we mean the sup norm of $u$ over the domain where $u$ is
defined. In this paper, we write $N = N(d,
\dots)$ if $N$ is a constant depending only on $d, \dots$.

\section{Main results}\label{main_sec}
                                                \label{section 3.13.1}

We consider the parabolic equation \eqref{para_main_eq}  with
coefficients $a^{jk}$, $b^{j}$, and $c$ satisfying the following
assumption.

\begin{assumption}\label{assum_01} 
The coefficients $a^{jk}$, $b^{j}$,
and $c$ are measurable functions defined on $\bR^{d+1}$, $a^{jk} =
a^{kj}$. There exist positive constants $\delta \in (0,1)$ and $K$ such
that
$$ |b^{j}(t, x)| \le K, \quad |c(t, x)| \le K,
$$
$$
\delta |\vartheta|^2 \le \sum_{j,k=1}^{d} a^{jk}(t,x) \vartheta^j
\vartheta^k
\le \delta^{-1} |\vartheta|^2
$$ for any $(t,x) \in \bR^{d+1}$ and $\vartheta \in \bR^d$.
\end{assumption}

We look for solutions of parabolic equations in the usual Sobolev space 
$$ 
W_p^{1,2}((S, T) \times \bR^d)=\{u:u,u_{t},u_{x},u_{xx}
\in L_{p}((S, T) \times \bR^d)\},
$$
$ -\infty \le S < T \le \infty$ with usual norm.
 Throughout the paper, as in \cite{Kr05}, we set 
$$
\Omega_T = (0,T) \times \bR^d.
$$ Thus, for instance,
$$ L_p(\Omega_T) = L_p((0,T) \times \bR^{d}),
\quad W_p^{1,2}(\Omega_T) = W_p^{1,2}((0,T) \times \bR^d).
$$

By $\WO{1,2}{p}(\Omega_T)$ we mean the collection of functions in
$W_p^{1,2}(\Omega_T)$ vanishing at $t = T$. We denote the differential
operator in \eqref{para_main_eq} by $L$, that is,
\begin{equation*} Lu = u_t + a^{jk} u_{x^j x^k} + b^{j} u_{x^j} + c u.
\end{equation*}

Our first result is about the case $p = 2$. In this case, we do not
require any regularity assumptions on the coefficients $a^{jk}$ if they
are functions of only $(t,x^1) \in \bR^{2}$.

\begin{theorem}                                   \label{para_main_02} 
Let Assumption \ref{assum_01} be
satisfied and let the coefficients $a^{ij}$ be
independent   of $x' \in \bR^{d-1}$.   Then for any $f \in
L_2(\Omega_T)$, there exists a unique $u \in \WO{1,2}{2}(\Omega_T)$
satisfying $Lu = f$ in $(0,T) \times \bR^d$. In addition, there is a
constant $N = N(d, \delta, K, T)$ such that, for any $u \in
\WO{1,2}{2}(\Omega_T)$,
$$
\| u \|_{W_{2}^{1,2}(\Omega_T)}
\le N \| Lu \|_{L_2(\Omega_T)}.
$$
\end{theorem}

\begin{remark} 
The assertion of Theorem \ref{para_main_02} is also valid  if 
 $a^{jk}(t,x)$  are uniformly continuous as functions  of $x' \in \bR^{d-1}$
uniformly in $(t,x^1) \in \bR^2$.  This can be shown by
using the standard techniques based on partitions of unity
and considering the equation on small time intervals
allowing one to absorb the $L_{2}$-norm into the $W^{1,2}_{2}$-norm.
Actually, there also is a  standard way, which can be found,
for instance,  in
\cite{Kr05}, to avoid
solving the equation step by step on small time intervals
moving down from $t=T$ to $t=0$.
\end{remark}

If $p \in (2,\infty)$, we suppose that the coefficients $a^{jk}$ are
measurable in $x^1 \in \bR$ and VMO in $(t,x') \in \bR^{d}$. To state
this assumption precisely, we introduce the following notation. Let 
$$ B_r(x) = \{ y \in \bR^d: |x- y| < r \},
$$
$$ B'_r(x') = \{ y' \in \bR^{d-1}: |x'- y'| < r \},
$$
$$ Q_r(t,x) = (t, t+r^2) \times B_r(x),
\quad
\Gamma_r(t,x') = (t, t+r^2) \times B'_r(x'),
$$
$$
\Lambda_{r}(t,x)= (t, t+r^2) \times (x^1-r, x^1+r) \times B'_r(x').
$$ Set $B_r = B_r(0)$, $B'_r = B'_r(0)$, $Q_r = Q_r(0)$, and so on.  By
$|B'_r|$ we mean  the $d-1$-dimensional volume of $B'_r(0)$. Denote $a =
(a^{jk})$ and 
$$
\text{osc}_{(t,x')} \left( a, \Lambda_r(t,x) \right)  = r^{-5}
|B'_r|^{-2}  \int_{x^1-r}^{x^1+r}  A_{(t,x')}(\tau) \, d\tau, 
$$ where
$$ A_{(t,x')}(\tau) = \int_{ (\sigma,y'), (\varrho,z') \in \Gamma_r(t,x')}
|a(\sigma, \tau, y') - a(\varrho, \tau, z') | \, dy' \, dz' \, d \sigma
\, d \varrho.
$$ Also denote
$$ a_R^{\#} = \sup_{(t,x) \in \bR^{d+1}} \sup_{r \le R} \,\,\,
\text{osc}_{(t,x')} \left( a, \Lambda_r(t,x) \right).
$$

\begin{assumption}\label{assum_02} There is a continuous function
$\omega(t)$ defined on $[0,\infty)$ such that $\omega(0) = 0$ and
$a_R^{\#} \le \omega(R)$ for all $R \in [0,\infty)$.
\end{assumption}

\begin{theorem}                                   \label{para_main_01}
Let $p \in (2, \infty)$ and let   Assumptions \ref{assum_01}  and
  \ref{assum_02} be satisfied. Then for any $f
\in L_p(\Omega_T)$, there exists a unique $u \in \WO{1,2}{p}(\Omega_T)$
such that $Lu = f$ in $(0,T) \times \bR^d$. Furthermore, there is a
constant $N = N(d, \delta, K, p, \omega, T)$ such that, for any $u \in
\WO{1,2}{p}(\Omega_T)$,
$$
\| u \|_{W_p^{1,2}(\Omega_T)}
\le N \| Lu \|_{L_p(\Omega_T)}.
$$
\end{theorem}

\begin{remark}
                                           \label{remark 7.23.1}
As usual in such situations,
from our proofs one can see that instead of the assumption
that
$a_R^{\#}\to0$ as $R\downarrow0$, actually, we are using that
there   exists $R\in(0,\infty)$ such that
$a_R^{\#}\leq\varepsilon$, where $\varepsilon>0$
is a constant depending only on other parameters
of the problem.
\end{remark}

We now show how to treat
 the Dirichlet and oblique derivative problems for parabolic equations in
half spaces. By the fact that coefficients are allowed to be measurable
in one direction,  in solving these problems, we need only the results
for equations in the whole space. Denote
$$ \bR^{d}_{+} = \{ x \in \bR^d : x^1 > 0 \},
\quad
\Omega^{+}_{T} = (0,T) \times \bR^{d}_{+},
$$
$$
\partial_t \Omega^{+}_{T}
= \{ (T,x) : x \in \bR^d_{+} \},
\quad
\partial_x \Omega^{+}_{T}
= \{ (t,0,x') : 0 \le t \le T, x' \in \bR^{d-1} \},
$$
$$
\partial' \Omega^{+}_{T} = \partial_t \Omega^{+}_{T} \cup  \partial_x \Omega^{+}_{T}.
$$

Below in this section
 we suppose that coefficients $a^{jk}$, $b^{j}$, and 
$c$ satisfy Assumption \ref{assum_01}.

\begin{theorem}                 
                                                          \label{theorem 3.29.1}
Let $2 \le p < \infty$. 
Assume that $a^{jk}$ are independent of $x' \in \bR^{d-1}$ if $p = 2$.
In case $p > 2$, we assume that $a^{jk}$ satisfy Assumption \ref{assum_02}.
Then for any $f \in L_p(\Omega^{+}_T)$, there exists a unique $u \in W_p^{1,2}(\Omega^{+}_{T})$ such that
$Lu = f$ in $(0,T) \times \bR^{d}_{+}$
and 
$u = 0$ on $\partial' \Omega^{+}_{T}$.
\end{theorem}

\begin{proof}
Introduce a new operator $\hat{L}v = \hat{a}^{jk} v_{x^j x^k} + \hat{b} v_{x^j} + \hat{c} v$, where $\hat{a}^{jk}$, $\hat{b}^{j}$, and $\hat{c}$ are defined as either even or odd extensions of $a^{jk}$, $b^{j}$, and $c$.
Specifically, 
for $j = k = 1$ or $j,k \in \{ 2, \dots, d \}$, even extensions:
$$
\hat{a}^{jk} = a^{jk}(t,x^1,x') 
\quad x^1 \ge 0,
\qquad
\hat{a}^{jk} = a^{jk}(t,-x^1,x')
\quad x^1 < 0.
$$
For $j = 2, \dots, d$, odd extensions:
$$
\hat{a}^{1j} = a^{1j}(t,x^1,x') 
\quad x^1 \ge 0,
\qquad
\hat{a}^{1j}= - a^{1j}(t,-x^1,x')
\quad x^1 < 0.
$$ Also set $\hat{a}^{j1} = \hat{a}^{1j}$. Similarly, $\hat{b}^1$ is the
odd extension of $b^1$, and $\hat{b}^{j}$, $j = 2, \dots, d$, and
$\hat{c}$ are even extensions of $b^{j}$ and $c$ respectively. We see
that the coefficients $\hat{a}^{jk}$, $\hat{b}^{j}$, and $\hat{c}$
satisfy Assumption \ref{assum_01}.  In addition, if $p>2$, the
coefficients $\hat{a}^{jk}$ satisfy Assumption \ref{assum_02} with $2
\omega$.

Let $\hat{f}$ be the odd extension of $f$.  Then it follows by Theorem
\ref{para_main_02} or \ref{para_main_01} that there exists a unique $u
\in \WO{1,2}{p}(\Omega_T)$ such that $\hat{L}u = \hat{f}$. It is easy to
check that $-u(t,-x^1,x') \in \WO{1,2}{p}(\Omega_T)$ also satisfy the
same equation, so by uniqueness we have $u(t,x^1,x') = -u(t,-x^1,x')$.
This and the fact that $u \in \WO{1,2}{p}(\Omega_T)$ show that that $u$,
as a function defined on $(0,T) \times \bR^{d}_{+}$, is a solution to $Lu
= f$ satisfying $u = 0$ on $\partial' \Omega^{+}_{T}$.

Uniqueness follows from the fact that the odd extension of a solution $u$
belongs to $\WO{1,2}{p}(\Omega_T)$ and the uniqueness of solutions to
equations in~$\Omega_T$.
\end{proof}

The following theorem addresses oblique derivative problems.

\begin{theorem}
Let $p$ and 
$a^{jk}$ be as in Theorem \ref{theorem 3.29.1}.
Let $\ell = (\ell^1, \cdots, \ell^d)$ be a vector in $\bR^d$ with $\ell^1 > 0$.
Then for any $f \in L_p(\Omega^{+}_T)$, there exists a unique $u \in W_p^{1,2}(\Omega^{+}_{T})$ satisfying
$Lu = f$ in $(0,T) \times \bR^{d}_{+}$,
$\ell^j u_{x^j} = 0$ on $\partial_x \Omega^{+}_{T} = 0$,
and
$u = 0$ on $\partial_t \Omega^{+}_{T} = 0$.
\end{theorem}

\begin{proof}
Let
$\varphi(x) = (\ell^1 x^1, \ell'x^1 + x')$, 
where $\ell' = (\ell^2, \dots, \ell^d)$.
Using this linear transformation and its inverse,
we reduce the above problem to a problem with Neumann boundary condition on $\partial_{x} \Omega^{+}_{T}$.
Note that, in case $p > 2$, the coefficients of the transformed equation satisfy Assumption \ref{assum_02} with $N \omega(N \cdot)$, where $N$ depends only on $d$ and $\ell$.
Then the latter problem is solved as in the proof of Theorem \ref{theorem 3.29.1} with the even extension of $f$.
\end{proof}

\begin{remark}
Solutions to problems in the above two theorems satisfy the $L_p$-estimate.
That is, if $u$ is a solution, then
$$
\| u \|_{W_p^{1,2}(\Omega^{+}_{T})} \le N \|f\|_{L_p(\Omega_{T}^{+})},
$$
where $N$ is a constant depending only on some or all 
 parameters -- $d$, $\delta$, $K$, $p$, $\omega$, $T$, $\ell$.
\end{remark}

\section{Proof of Theorem \ref{para_main_02}} 
                                                         \label{L_2_case}

Introduce
\begin{equation}
                                                   \label{Def_L_0} 
L_0 u(t,x) = u_t(t,x) + a^{jk}(t,x^1)
u_{x^j x^k}(t,x),
\end{equation}

\begin{lemma}
                                                      \label{lemma 3.10.1} 
Assume that $d=1$. Then   for any $\lambda>0$ and $f\in L_{2}(\bR^{2})$
there exists a unique solution $u\in W^{1,2}_{2}(\bR^{2})$ of 
the equation
$L_{0}u-\lambda u=f$. Furthermore, there is a constant $N=N(\delta)$
such that for any $\lambda\geq0$ and
$u\in W^{1,2}_{2}(\bR^{2})$ we have
$$
\|u_{t}\|_{L_{2}(\bR^{2})}+\|u_{xx}\|_{L_{2}(\bR^{2})}
+\sqrt{\lambda}\|u_{x }\|_{L_{2}(\bR^{2})}
$$
\begin{equation}
                                                         \label{3.10.1}
+\lambda\|u \|_{L_{2}(\bR^{2})}\leq
N\|L_{0}u-\lambda u\|_{L_{2}(\bR^{2})}.
\end{equation}
\end{lemma}

\begin{proof} 
As usual we only need prove \eqref{3.10.1}
and only for $u\in C^{\infty}_{0}(\bR^{2})$. Take such a function,
denote $a=a^{11}$, $f:=L_{0}u-\lambda u$, and write
$$
a^{-1/2}f=a^{1/2}u_{xx}+a^{-1/2}(u_{t}-\lambda u),
$$
$$
 a^{-1}f^{2}=au_{xx}^{2}+2u_{xx}(u_{t}-\lambda u)
+a^{-1}(u_{t}-\lambda u)^{2}.
$$
Then integrate through the last equation over $\bR^{2}$
and notice that
$$
2\int_{\bR^{2}}u_{xx} u_{t}\,dxdt=
-2\int_{\bR^{2}}u_{x } u_{xt}\,dxdt=
-\int_{\bR^{2}}\frac{\partial}{\partial t}u_{x }^{2}\,dtdx=0,
$$
$$
2\int_{\bR^{2}}u_{xx} u \,dxdt=
-\int_{\bR^{2}} u_{x }^{2}\,dxdt.
$$
Then we find
$$
\delta^{-1}\int_{\bR^{2}}f^{2}\,dxdt\geq
\delta\int_{\bR^{2}}u_{xx}^{2}\,dxdt+
\lambda\int_{\bR^{2}}u_{x}^{2}\,dxdt+
\delta\int_{\bR^{2}}(u_{t}-\lambda u)^{2}\,dxdt.
$$
Upon observing that
$$
\int_{\bR^{2}}(u_{t}-\lambda u)^{2}\,dxdt
=\int_{\bR^{2} }u_{t} ^{2}\,dxdt
-2\lambda\int_{\bR^{2}} u_{t}  u \,dxdt
+\lambda^{2}\int_{\bR^{2}}u ^{2}\,dxdt,
$$
$$
 2 \int_{\bR^{2}} u_{t}  u \,dxdt
=-\int_{\bR^{2}}\frac{\partial}{\partial t}u ^{2}\,dxdt=0
$$
we finish the proof.
\end{proof}

We now generalize Lemma \ref{lemma 3.10.1}
to cover the multidimensional case.
\begin{theorem}
                                                  \label{theorem 3.10.2} 
 For any $\lambda>0$ and $f\in L_{2}(\bR^{d+1})$
there exists a unique solution $u\in W^{1,2}_{2}(\bR^{d+1})$ of 
the equation
$L_{0}u-\lambda u=f$. Furthermore, there is a constant $N=N( \delta)$
such that for any $\lambda\geq0$ and
$u\in W^{1,2}_{2}(\bR^{d+1})$ we have
$$
\|u_{t}\|_{L_{2}(\bR^{d+1})}+\|u_{xx}\|_{L_{2}(\bR^{d+1})}
+\sqrt{\lambda}\|u_{x }\|_{L_{2}(\bR^{d+1})}
$$
\begin{equation}
                                                         \label{3.10.2}
+\lambda\|u \|_{L_{2}(\bR^{d+1})}\leq
N\|L_{0}u-\lambda u\|_{L_{2}(\bR^{d+1})}.
\end{equation}
\end{theorem}

It is worth saying that by 
$$
\|u_{x}\|_{L_{2}(\bR^{d+1})}\quad\text{and}\quad
\|u_{xx}\|_{L_{2}(\bR^{d+1})}
$$
in \eqref{3.10.2} we mean $L_{2}$-norms of
$$
 \big(\sum_{k}|u_{x^{k}}|^{2}\big)^{1/2} \quad\text{and}\quad
\big(\sum_{k,j}|u_{x^{k}x^{j}}|^{2}\big)^{1/2},
$$
respectively. Different definitions could make $N$ depend also on $d$.

We prove this theorem after some preparations.
Again it suffices to only prove \eqref{3.10.2} and only for 
$u \in C^{\infty}_{0}(\bR^{d+1})$. In addition
we may assume that $a^{ij}$ are infinitely differentiable.
Fix such $u$, $a^{ij}$, and $\lambda\geq0$ and set 
$$
f:=
L_{0}u-\lambda u.
$$

 Let $\xi \in \bR^{d-1}$ and let $\tilde{\psi}(t,x^1,\xi)$ denote
 the Fourier
transform of $\psi(t,x^1,x')$ with respect to $x' \in \bR^{d-1}$. By
taking the Fourier transform (with respect to $x' \in \bR^{d-1}$),
 we obtain
$$
\tilde{u}_t(t,x^1,\xi) + \coa(t,x^1) \tilde{u}_{x^1x^1}(t,x^1,\xi) +
\mathrm{i} \, 2 \, \cob(t,x^1,\xi) \tilde{u}_{x^1}(t,x^1,\xi)
$$
\begin{equation}
                                                        \label{3.11.1}
 - \coc(t,x^1,\xi) \tilde{u}(t,x^1,\xi) -
\lambda \tilde{u}(t,x^1,\xi) =\tilde{f}(t,x^1,\xi),
\end{equation}
where $\mathrm{i} = \sqrt{-1}$,
$$
\coa(x^1) = a^{11}(x^1), 
\qquad
\cob(t,x^1,\xi) = \sum_{j=2}^d a^{1j}(t,x^1) \xi^j,
$$
$$
\coc(t,x^1,\xi) = \sum_{j,k=2}^d a^{jk}(t,x^1) \xi^j \xi^k.
$$ 
 Introduce a function 
$$
\rho(t,x^1,\xi) = \tilde{u}(t,x^1,\xi) \, e^{\mathrm{i} \phi(t,x^1,\xi)},
$$
where $\phi(t,0, \xi) = 0$ and $\phi_{x^1}(t,x^1, \xi) =
 \coa^{-1}\cob (t,x^1,\xi)$. It is easy to see that $\rho$ satisfies
\begin{equation}
                                                                                            \label{3.10.6}
 \rho_t + \coa\rho_{x^1x^1}-\left(\coc- \coa^{-1}\cob ^2  
+ \lambda + \mathrm{i}  \phi_t + \mathrm{i}\coa
\phi_{x^1x^1} \right) \rho = \tilde{f} e^{\mathrm{i} \phi}.
\end{equation}

In the following lemma $\xi$ is considered
 as a parameter.

\begin{lemma}
                                               \label{lemma 3.10.3}
 Let $|\xi|^2 + \lambda > 0$. 
Then we have
\begin{equation}
                                                                                            \label{3.10.3}
|\rho(t,x^1,\xi)|\leq\hat{\rho}(t,x^1,\xi),
\end{equation}
where, for each $\xi\in\bR^{d-1}$, $\hat{\rho}(t,x^1,\xi)$ 
is the unique $W^{1,2}_{2}(\bR^{2})$
solution of 
 \begin{equation}
                                      \label{3.10.4}
\hat{\rho}_t + \coa\hat{\rho}_{x^1x^1} -( \lambda
 + \delta^{3}|\xi|^{2} )\hat{\rho} = -|\tilde{f}|  .
\end{equation}

In particular, (by Lemma \ref{lemma 3.10.1})
$$
( |\xi|^{2}
+\lambda) \|\tilde{u}(\cdot,\cdot,\xi)\| _{L_{2}(\bR^{2})}=
( |\xi|^{2}
+\lambda) \|\rho(\cdot,\cdot,\xi)\| _{L_{2}(\bR^{2})}
$$
 \begin{equation}                             \label{3.10.10}
\leq
( |\xi|^{2}
+\lambda) \|\hat\rho(\cdot,\cdot,\xi)\| _{L_{2}(\bR^{2})}
\leq N(\delta)\|\tilde{f}(\cdot,\cdot,\xi)\| _{L_{2}(\bR^{2})}.
\end{equation}
\end{lemma}

\begin{proof}
First, observe that by Lemma \ref{lemma 3.10.1}
the function $\hat{\rho}$ indeed exists and by the maximum principle
it is nonnegative. Also since $|\tilde{f}|$
is  Lipschitz continuous, $\hat{\rho}$
is twice continuously differentiable in $x$
and once in $t$.

Assume that, for a fixed $\xi$, \eqref{3.10.3} is violated.
Then, due to the fact that
 $\rho$ has a compact support, there is a point
$(t_{0},x^{1}_{0})$ such that
 \begin{equation}
                                          \label{3.10.5} 
| \rho
(t_{0},x^{1}_{0})|-
\hat{\rho}(t_{0},x^{1}_{0})=\max_{\bR^{2}}
(| \rho(t,x^{1})|-
\hat{\rho}(t,x^{1}))>0.
\end{equation}
Since $| \rho (t_{0},x^{1}_{0})|>0$ and $\rho$ is smooth,
 the function $|\rho|$ is twice differentiable at $(t_{0},x_{0})$
 and at this point
 $$
|\rho|_{x^{1}}=\frac{\Re(\bar{\rho}\rho_{x^{1}})}
{|\rho|}=\hat{\rho}_{x^{1}},\quad 
|\rho|_{t}=\frac{\Re(\bar{\rho}\rho_{t})}
{|\rho|}=\hat{\rho}_{t},
$$
$$
|\rho|_{x^{1}x^{1}}=\frac{1}{|\rho|^{3}}
\big(|\rho|^{2}|\rho_{x^{1}}|^{2}-
(\Re(\bar{\rho}\rho_{x^{1}}))^{2}\big)+
\frac{1}{|\rho| }\Re(\bar{\rho}\rho_{x^{1}x^{1}})
\leq\hat{\rho}_{x^{1}x^{1}}.
$$
Obviously, $(\Re(\bar{\rho}\rho_{x^{1}}))^{2}
\leq|\rho|^{2}|\rho_{x^{1}}|^{2}$, so that we also have
$$
\frac{1}{|\rho| }\Re(\bar{\rho}\rho_{x^{1}x^{1}})
\leq\hat{\rho}_{x^{1}x^{1}}.
$$

Next, we multiply 
\eqref{3.10.6} by $\eta:=\bar{\rho}/|\rho|$ and
take real parts of both sides
to get
\begin{equation}
                                                        \label{3.10.7}
 \Re(\eta\tilde{f} e^{\mathrm{i} \phi})=\Re(\eta\rho_t) + \coa
\Re(\eta\rho_{x^1x^1} )- (\coc -  \coa^{-1}\cob ^2 
+\lambda )|\rho|  
\end{equation}

Concentrate on this equation
at the point $(t_{0},x_{0})$ 
 and use the above manipulations with the derivatives
to see that  at $(t_{0},x_{0})$
$$
\Re(\eta\tilde{f} e^{\mathrm{i} \phi})
\leq\hat{\rho}_{t}+\coa\hat{\rho}_{x^{1}x^{1}}
- (\coc -  \coa^{-1}\cob ^2 
+\lambda )|\rho|
$$
$$
=-|\tilde{f}|+( \delta^{3} |\xi|^2   +\lambda)\hat{\rho}
- (\coc -  \coa^{-1}\cob ^2 
+\lambda )|\rho|.
$$
Here, $|\rho|>\hat{\rho}\geq0$ and as is easy to check (see, for instance,
Lemma 3.1 in ~\cite{KiK}),
$\coc -  \coa^{-1}\cob ^2 \geq  \delta^{3}|\xi|^2$. Therefore,
(always at $(t_{0},x_{0})$)
$$
(\coc -  \coa^{-1}\cob ^2 
+\lambda )|\rho|>(\delta^{3}|\xi|^2+\lambda)\hat{\rho},
$$
so that we get
$$
\Re(\eta\tilde{f} e^{\mathrm{i} \phi})\leq
 -|\tilde{f}|+( \delta^{3}|\xi|^2 + \lambda)\hat{\rho}
- (\coc -  \coa^{-1}\cob ^2 
+\lambda )|\rho|<-|\tilde{f}|.
$$
This leads to a contradiction because $|\eta|=1$
and proves the lemma.\end{proof}

\begin{lemma}
                                               \label{lemma 3.10.4}
 For any $\varepsilon>0$,
there exists a constant $ N(\varepsilon,\delta)$ such that
 \begin{equation}
                               \label{3.10.8}
 (|\xi|+ \sqrt{\lambda})
 \|\rho_{x}(\cdot,\cdot,\xi)\| _{L_{2}(\bR^{2})}
\leq N(\varepsilon,\delta)
\|\tilde{f}(\cdot,\cdot,\xi)\| _{L_{2}(\bR^{2})}+\varepsilon
\|\tilde{u}_{t}(\cdot,\cdot,\xi)\| _{L_{2}(\bR^{2})}.
\end{equation}
\end{lemma}

\begin{proof}

We go back to equation \eqref{3.10.7}, which we multiply
by $ |\rho|$, divide by $\coa$, 
then integrate over $\bR^{2}$,
and use that $\coc\leq \delta^{-1}|\xi|^{2}$
  and $|b| \le \delta^{-1}|\xi|$. We also use the fact that
$$
\Re\int_{\bR^{2}}\bar{\rho}\rho_{x^{1}x^{1}}\,dxdt
= -  \int_{\bR^{2}}|\rho_{x^{1}}|^{2}\,dxdt,
$$
$$
2\Re(\bar{\rho}\rho_{t})=\frac{\partial}{\partial t}
|\rho|^{2}=\frac{\partial}{\partial t}
|\tilde{u}|^{2}=2\Re(\bar{\tilde{u}}\tilde{u}_{t}).
$$
Then we obtain
$$
\int_{\bR^{2}}|\rho_{x^{1}}|^{2}\,dxdt
\leq2\int_{\bR^{2}}\coa^{-1}|\bar{\tilde{u}}\tilde{u}_{t}|\,dxdt
$$
$$
+N(\lambda+|\xi|^{2})\int_{\bR^{2}}|\rho|^{2}\,dxdt
+2\int_{\bR^{2}} \coa^{-1} |\rho\tilde{f}| \,dxdt.
$$
We estimate the terms on the right by using Young's inequality
and assuming without losing generality that
$\lambda+|\xi|^{2}\ne0$.
For instance,
$$
2\int_{\bR^{2}}\coa^{-1}|\bar{\tilde{u}}\tilde{u}_{t}|\,dxdt
\leq  \delta^{-2} \varepsilon^{-1}
(\lambda+|\xi|^{2}) \int_{\bR^{2}}|\tilde{u}|^{2}\,dxdt
$$
$$
+\varepsilon(\lambda+|\xi|^{2})^{-1}
\int_{\bR^{2}}|\tilde{u}_{t}|^{2}\,dxdt.
$$
We also use \eqref{3.10.10}. Then we easily get
\eqref{3.10.8}.\end{proof}

\begin{proof}[{\bf Proof of Theorem \ref{theorem 3.10.2}}]
Since 
$$
\tilde{u}=\rho e^{- \mathrm{i} \phi} ,\quad
 \tilde{u}_{x^{1}}=[\rho_{x^{1}}-\mathrm{i}
\coa^{-1}\cob\rho]e^{-\mathrm{i}\phi},
$$
and $|\cob|\leq N|\xi|$, Lemmas \ref{lemma 3.10.3}
 and \ref{lemma 3.10.4} imply that for any
$\varepsilon>0$ there is an $ N(\varepsilon,\delta)$ such that
$$
(|\xi|^{4}+\lambda^{2})\|\tilde{u}(\cdot,\cdot,\xi)\|^{2}
_{L_{2}(\bR^{2})}
+(|\xi|^{2}+\lambda)\|\tilde{u}_{x^{1}}(\cdot,\cdot,\xi)\|^{2}
_{L_{2}(\bR^{2})}
$$
\begin{equation}
                                                     \label{3.11.3}
\leq N(\varepsilon,\delta)
\|\tilde{f}(\cdot,\cdot,\xi)\|^{2}
_{L_{2}(\bR^{2})}+\varepsilon\|\tilde{u}_{t}(\cdot,\cdot,\xi)\|^{2}
_{L_{2}(\bR^{2})}.
\end{equation}
Then \eqref{3.11.1} shows that for  any
$\varepsilon>0$ there is an $ N(\varepsilon,\delta)$ such that
$$
\|(\tilde{u}_{t}+\coa\tilde{u}_{x^{1}x^{1}})(\cdot,\cdot,\xi)\|^{2}
_{L_{2}(\bR^{2})}\leq N(\varepsilon,\delta)
\|\tilde{f}(\cdot,\cdot,\xi)\|^{2}
_{L_{2}(\bR^{2})}+\varepsilon\|\tilde{u}_{t}(\cdot,\cdot,\xi)\|^{2}
_{L_{2}(\bR^{2})},
$$
which after being combined with Lemma \ref{lemma 3.10.1}
(with $\lambda=0$ there)
leads to
$$
\|\tilde{u}_{t} (\cdot,\cdot,\xi)\|^{2}
_{L_{2}(\bR^{2})}+
\| \tilde{u}_{x^{1}x^{1}}(\cdot,\cdot,\xi)\|^{2}
_{L_{2}(\bR^{2})}
$$
\begin{equation}
                                             \label{3.21.1}
\leq N(\varepsilon,\delta)
\|\tilde{f}(\cdot,\cdot,\xi)\|^{2}
_{L_{2}(\bR^{2})}+\varepsilon\|\tilde{u}_{t}(\cdot,\cdot,\xi)\|^{2}
_{L_{2}(\bR^{2})}.
\end{equation}
The reader might have noticed that in the above computations
the constants $N(\varepsilon,\delta)$ are changing from line to line
and $\varepsilon$ was sometimes multiplied by a constant
of type $N(\delta)$. However, $N(\delta)\varepsilon$ is as arbitrary
as $\varepsilon$. Upon taking $\varepsilon=1/2$ in 
\eqref{3.21.1}
we conclude that
\begin{equation}
                                                     \label{3.11.2}
\|\tilde{u}_{t} (\cdot,\cdot,\xi)\|^{2}
_{L_{2}(\bR^{2})}+
\| \tilde{u}_{x^{1}x^{1}}(\cdot,\cdot,\xi)\|^{2}
_{L_{2}(\bR^{2})}
\leq N 
\|\tilde{f}(\cdot,\cdot,\xi)\|^{2}
_{L_{2}(\bR^{2})}.
\end{equation}

After that \eqref{3.11.3} yields
$$
(|\xi|^{4}+\lambda^{2})\|\tilde{u}(\cdot,\cdot,\xi)\|^{2}
_{L_{2}(\bR^{2})}
+(|\xi|^{2}+\lambda)\|\tilde{u}_{x^{1}}(\cdot,\cdot,\xi)\|^{2}
_{L_{2}(\bR^{2})}
$$
\begin{equation}
                                                     \label{3.11.4}
\leq
N 
\|\tilde{f}(\cdot,\cdot,\xi)\|^{2}
_{L_{2}(\bR^{2})}.
\end{equation}

To get \eqref{3.10.2}  now it only remains to
integrate through \eqref{3.11.2} and \eqref{3.11.4}
 with respect to $\xi$ and use Parseval's identity.
The theorem is proved.
\end{proof}

 Theorem~\ref{para_main_02} is derived from Theorem 
\ref{theorem 3.10.2} in a  standard way, which can be found,
for instance,  in
\cite{Kr05}. Theorem ~\ref{para_main_02} is proved.

\section{Auxiliary results for equation in $W_p^{1,2}(\bR^{d+1})$}
\label{aux_W_p}

We assume in this section that
$a^{jk}$ are measurable functions only of  $x^1 \in \bR$. Set
$$ L_0 u(t,x) = u_t(t,x) + a^{jk}(x^1) u_{x^j x^k}(t,x).
$$ By $\partial' Q_r(t,x)$ we mean the parabolic boundary of $Q_r(t,x)$
defined as
$$
\partial' Q_r(t,x) = \left([t,t+r^2] \times \partial B_r(x)\right)
\cup \{ (t+r^2,y) : y \in B_r(x) \}.
$$

\begin{lemma} \label{para_L_2_esti} There exists $N = N(d,\delta)$ such
that, for $u \in W_2^{1,2}(Q_r)$ with $u|_{\partial' Q_r} = 0$, we have
\begin{equation}				\label{Q_R_L_2_esti} r^2 \int_{Q_r} |u_x|^2 \, dx \,
dt + \int_{Q_r} |u|^2 \, dx \, dt
\le N r^4 \int_{Q_r} |L_0 u|^2 \, dx \, dt.
\end{equation}
\end{lemma}

\begin{proof} Assume that \eqref{Q_R_L_2_esti} is true when $r = 1$. For
$u \in W_2^{1,2}(Q_r)$ with $u|_{\partial' Q_r} = 0$, we set
$$
\hat{L}_0 = \frac{\partial}{\partial t} + a^{jk}(r x^1)
\frac{\partial^2}{\partial x^j \partial x^k}
\quad
\text{and}
\quad
\hat{u}(t,x) = r^{-2} u(r^2 t, r x).
$$ Then $\hat{u} \in W_2^{1,2}(Q_1)$ and $\hat{L}_0 \hat{u}(t,x)
 = L_{0}
u(r^2 t, r x)$ in $Q_1$. Since $\hat{L}_0$ satisfies the same ellipticity
condition as $L_0$ does, we have
$$
\int_{Q_r} |u|^2 \, dx \, dt = r^{d+6} \int_{Q_1} |\hat{u}|^2 \, dx \, dt
$$
$$
\le N r^{d+6} \int_{Q_1} |\hat{L}_0 \hat{u}|^2 \, dx \, dt = N r^4
\int_{Q_r} |L_0 u|^2 \, dx \, dt.
$$ Also
$$
\int_{Q_r} |u_x|^2 \, dx \, dt = r^{d+4} \int_{Q_1} |\hat{u}_x|^2 \, dx
\, dt
$$
$$
\le N r^{d+4} \int_{Q_1} |\hat{L}_0 \hat{u}|^2 \, dx \, dt = N r^2
\int_{Q_r} |L_0 u|^2 \, dx \, dt.
$$ This shows that we need to prove the lemma only for $r = 1$.

In this case, we divide $L_0$ by $a^{11}(x^1)$. That is, by setting 
$$
f :=
u_t + a^{jk} u_{x^j x^k},\quad \hat{a}^{jk} := a^{jk}/a^{11},
$$
 we have
$$
u_t/a^{11} + \hat{a}^{jk} u_{x^j x^k} = f/a^{11}.
$$
 Then using the
ellipticity of $a^{jk}$ and integration by parts, we obtain
$$
\delta^2 \int_{Q_1} |u_x|^2 \, dx \, dt
\le \int_{Q_1} \hat{a}^{jk} u_{x^j} u_{x^k} \, dx \, dt = - \int_{Q_1} u
\, \hat{a}^{jk} u_{x^j x^k} \, dx \, dt
$$
$$ = \int_{Q_1} \frac{u}{a^{11}} \left( u_t - f \right) \, dx \, dt.
$$ Note that
$$
\int_{Q_1} \frac{u}{a^{11}} u_t \, dx \, dt = \int_{B_1} \frac{1}{a^{11}}
\int_0^1 u \, u_t \, dt \, dx = - \int_{B_1} \frac{1}{2 a^{11}} u(0,x)^2
\, dx \le 0,
$$ where we used the fact that $a^{11}$ is independent of $t$ and $u(1,x)
= 0$. Thus we have
$$
\delta^2 \int_{Q_1} |u_x|^2 \, dx \, dt
\le - \int_{Q_1} \frac{u}{a^{11}} f \, dx \, dt
$$
$$
\le \delta^{-1} \left( \int_{Q_1} u^2 \, dx \, dt \right)^{1/2}
\left( \int_{Q_1} f^2 \, dx \, dt \right)^{1/2}.
$$ By using Poincar\'{e}'s inequality, we estimate the integral of $u^2$
in the last term through that of $|u_x|^2$.  This gives us the needed
estimate for $u_x$. For the estimate for $u$, we use Poincar\'{e}'s
inequality once again. The lemma is proved.
\end{proof}

\begin{lemma} \label{glo_loc} Let $0 < r < R$. There exists $N =
N(d,\delta)$ such that, for $u \in W_{2}^{1,2}(Q_R)$,
$$
\|u\|_{W_2^{1,2}(Q_r)}
\le N \left( \| L_0 u - u \|_{L_2(Q_R)} + (R-r)^{-2} \| u \|_{L_2(Q_R)}
\right).
$$
\end{lemma}

\begin{proof} The proof is just a repetition
of the proof  of Lemma 4.2 in
\cite{KiK}
and is based on 
 \eqref{3.10.2}
with  $\lambda=1$ and $Q_{m}$ and $\zeta_{m}$, specified below.
Introduce
$$ Q_m = Q_{r_m} = (0, {r_m}^2) \times B_{r_m},
\quad m = 1, 2, \dots,
$$ where
$r_0 = r$ and $r_m = r + (R-r) \sum_{k=1}^{m} 2^{-k}$. Also let $\zeta_m
\in C_0^{\infty}(\bR^{d+1})$ be such that
$$
\zeta_m(t,x) = 
\left\{
\begin{aligned} 1 \quad &\text{on} \quad Q_m \\  0 \quad &\text{on} \quad
\bR^{d+1} \setminus \left[(-{r_{m+1}}^2, {r_{m+1}}^2) \times
B_{r_{m+1}}\right]
\end{aligned}
\right.
$$ and
$$ | (\zeta_m)_x |_0 \le N \frac{2^{m+1}}{R-r},
\,\,\, | (\zeta_m)_t |_0 \le N \frac{2^{2m+2}}{(R-r)^2},
\,\,\, | (\zeta_m)_{xx} |_0 \le N \frac{2^{2m+2}}{(R-r)^2},
$$ where $N$ is a constant.  In fact, we construct $\zeta_m$ as follows.
Let $g(t)$ be an infinitely differentiable function defined on $\bR$ such
that $g(t) = 1$ for $t \le 1$, $g(t) = 0$, for $t \ge 2$, and $0 \le g
\le 1$.  Then define
$$
\rho_m(x) = g( 2^{m+1}(R-r)^{-1} (|x| - r_m) + 1 ),
$$
$$
\eta_m(t) = g( 2^{m+1}(R-r)^{-1} (\sqrt{|t|} - r_m) + 1 ),
$$
$$
\zeta_m(t,x) = \eta_m(t) \rho_m(x).
$$
\end{proof}

\begin{lemma} \label{pro_har} Let $0 < r < R$ and $\gamma = (\gamma^1,
\cdots, \gamma^d)$ be a multi-index such that $\gamma^1 = 0, 1, 2$.  If
$h$ is a sufficiently smooth function defined on $Q_R$ such that $L_0 h =
0$ in $Q_R$, then 
$$
\int_{Q_r} | D_{t}^m D_{x}^{\gamma} h |^2 \, dx \, dt
\le N \int_{Q_R} |h|^2 \, dx \, dt,
$$ where $m$ is a nonnegative integer and $N = N(d,\delta, \gamma, m, R,
r)$.
\end{lemma}

\begin{proof} Since $a^{jk}$ are independent of $t \in \bR$ and $x' \in
\bR^{d-1}$, we have $L_0 (D^m_t h) = 0$ and $L_0 (D^{\gamma'}_{x} h) =
0$,  where $\gamma' = (0, \gamma^2, \dots, \gamma^d)$. Then the proof is
completed using Lemma \ref{glo_loc}  and the argument in the proof of
Lemma 4.4 in \cite{KiK}.
\end{proof}

Throughout the rest of this paper, depending on the context,
 by $h_{x'}$ we mean 
one of
$h_{x^j}$, $j =
2, \dots, d$ or the whole collection consisting of them. 
By $h_{x}$ we mean one of
$h_{x^j}$, $j =
1, \dots, d$ or
the full gradient of $h$ with respect to $x$. Also,
by $h_{xx'}$ we mean one of  $h_{x^j x^k}$, where 
$j \in \{1, \dots, d\}$ and
$k \in \{2, \dots, d\}$ or
the collection of them. Norms of these collections
are defined arbitrarily.

\begin{lemma} \label{sup_holder_esti} Let $h$ be a sufficiently smooth
function defined on $Q_4$ such that $L_0 h = 0$ in $Q_4$.  Then
$$
\sup_{Q_1}|h_{tt}| + \sup_{Q_1}|h_{tx}| + \sup_{Q_1}|h_{tx'x}| +
\sup_{Q_1}|h_{x'xx}|
\le N \| h \|_{L_2(Q_3)},
$$ where $N = N(d, \delta)$.
\end{lemma}

\begin{proof} We prove that 
\begin{equation}\label{ess_sup}
\sup_{Q_1}|h| + \sup_{Q_1}|h_{x^1}|
 \le N \| h \|_{L_2(Q_{r})},
\end{equation} where $2 < r < 3$
and $N=N(r,d,\delta)$. 
If this is true, then using the fact
that 
$$
L_0 h_{t} =L_0 h_{tt} =L_0 h_{tx'} =L_0 h_{tx'x'} =
L_0 h_{x'x'} = L_0 h_{x'x'x'} =0
$$
 we obtain
$$
\sup_{Q_1}|h_{tt}| + \sup_{Q_1}|h_{tx}| + \sup_{Q_1}|h_{t x'x}| +
\sup_{Q_1}|h_{x'x'x}| 
\le N \sum_{k + |\gamma| \le 3} \| D^{k}_{t} D^{\gamma}_{x'} h
\|_{L_2(Q_r)}.
$$ 
This and Lemma \ref{pro_har} prove all the 
desired estimates except
$$
\sup_{Q_1}|h_{x'x^1x^1}| \le N \| h \|_{Q_3}.
$$ However, this one holds true as well because
$$ 
a^{11}h_{x'x^1x^1} = - h_{x't} - \sum_{j \ne 1 \, \text{or} \, k \ne 1}
a^{jk} h_{x'x^j x^k}.
$$

To prove \eqref{ess_sup}, we observe that,  due to the Sobolev embedding
theorem, there exist positive constants $m$ and $N$ such that
$$
\sup_{Q_1} |h_{x^1}|
\le N \sum_{k + |\gamma| \le m} 
\left( \| D_t^k D_{x'}^{\gamma} h_{x^1} \|_{L_2(Q_2)} + \| D_t^k
D_{x'}^{\gamma} h_{x^1 x^1} \|_{L_2(Q_2)} \right).
$$ By Lemma \ref{pro_har}, the right side of the above inequality is not
greater than a constant
$N=N(r,d,\delta)$
times $\| h \|_{L_2(Q_r)}$, $2 < r < 3$. This
proves that 
$$
\sup_{Q_1} |h_{x^1}| \le N \| h \|_{L_2(Q_r)}.
$$ Similarly, we have the same inequality as above with $h$ in place of
$h_{x^1}$. Therefore,  \eqref{ess_sup} is proved, so is the
lemma. ~\end{proof}

Let $u \in C_0^{\infty}(\bR^{d+1})$. Assume that $a^{jk}(x^1)$ are
infinitely differentiable.  Then there exists a sufficiently smooth
function $h$ defined on $Q_4$ such that
$$
\left\{ 
\begin{aligned} L_0 h &= 0 \quad \text{in} \quad Q_4 \\ h  &= u \quad
\text{on} \quad \partial' Q_4
\end{aligned}\right..
$$ The functions $u$ and $h$ satisfy the following inequality.

\begin{lemma} \label{sup_h_f} There exists a constant $N = N(d, \delta)$
such that
$$
\sup_{Q_1}|h_{tt}| + \sup_{Q_1}|h_{tx}| + \sup_{Q_1}|h_{tx'x}| +
\sup_{Q_1}|h_{x'xx}|
$$
$$
\le N  \left( \| L_0 u \|_{L_2(Q_4)} + \| u_{xx} \|_{L_2(Q_4)} \right).
$$
\end{lemma}

\begin{proof} We need only follow the argument in Lemma 4.6 in
\cite{KiK} along with Lemma
\ref{para_L_2_esti} and \ref{sup_holder_esti}.
\end{proof}

Denote by $(u)_{Q_r(t_0,x_0)}$ the average value of a function $u$ over
$Q_r(t_0,x_0)$, that is,
$$ (u)_{Q_r(t_0,x_0)} = \dashint_{Q_r(t_0,x_0)} u(t,x) \, dx \, dt.
$$

\begin{lemma} \label{h_final_esti} Let $\kappa \ge 4$ and $r > 0$. Assume
that $a^{jk}(x^1)$ are infinitely differentiable. For $u \in
C_0^{\infty}(\bR^{d+1})$, we find a smooth function $h$ defined on
$Q_{\kappa r}$ such that $L_0 h = 0$ in $Q_{\kappa r}$ and $h = u$ on
$\partial' Q_{\kappa r}$.  Then there exists a constant $N = N(d,
\delta)$ such that
$$
\dashint_{Q_r} | h_{t} - (h_{t})_{Q_r} |^2 \, dx \, dt + \dashint_{Q_r} |
h_{xx'} - (h_{xx'})_{Q_r} |^2 \, dx \, dt$$
$$
\le N \kappa^{-2} 
\left[ (|L_0 u|^2)_{Q_{\kappa r}} + (|u_{xx}|^2)_{Q_{\kappa r}} \right].
$$
\end{lemma}

\begin{proof} By the dilation argument as in the proof of Lemma
\ref{para_L_2_esti},  we need to prove our assertion
 only in the case $r = 1$. In
this case, we use Lemma \ref{sup_h_f} and the dilation argument again  to
obtain
\begin{multline}\label{sup_h_kappa}
\kappa^{2}\sup_{Q_{\kappa/4}}|h_{tt}|^2  + 
\sup_{Q_{\kappa/4}}|h_{tx}|^2  +\kappa^{2}
\sup_{Q_{\kappa/4}}|h_{tx'x}|^2 + \sup_{Q_{\kappa/4}}|h_{x'xx}|^2 \\
\le N \kappa^{-2} 
\left[ (|L_0 u|^2)_{Q_{\kappa}} + (|u_{xx}|^2)_{Q_{\kappa}} \right].
\end{multline} Set $v$ to be either $h_t$ or $h_{xx'}$. Then by the fact
that $\kappa \ge 4$ it follows that
$$
\dashint_{Q_1} | v - (v)_{Q_1} |^2 \, dx \, dt 
\le N \sup_{Q_{\kappa/4}}|v_t|^{2} + N \sup_{Q_{\kappa/4}}|v_x|^{2}.
$$
 This and \eqref{sup_h_kappa} prove the assertion 
of the lemma in
 case $r =1$. The lemma is
proved. ~\end{proof}

\begin{lemma} \label{u_sharp} There exists a constant $N = N(d, \delta)$
such that, for any $\kappa \ge 4$, $r > 0$, and $u \in
C_0^{\infty}(\bR^{d+1})$, we have
$$
\dashint_{Q_r} |u_{t} - (u_{t})_{Q_r}|^2 \, dx \, dt + \dashint_{Q_r} |u_{xx'}
- (u_{xx'})_{Q_r}|^2 \, dx \, dt
$$
$$
\le N \kappa^{d+2} \left( |L_0 u|^2 \right)_{Q_{\kappa r}} + N
\kappa^{-2} \left( |u_{xx}|^2 \right)_{Q_{\kappa r}}.
$$
\end{lemma}

\begin{proof} Use Lemma \ref{h_final_esti}, \ref{glo_loc},
\ref{para_L_2_esti}, and the argument in the proof of Lemma 4.8 in
\cite{KiK} (also see Remark 4.3
there).~\end{proof}

\section{Proof of Theorem \ref{para_main_01}} \label{proof_th_p}

We assume in this section that all assumptions 
of Theorem \ref{para_main_01}
are satisfied. 
However, in Theorem \ref{main_sharp}
 the assumption that $\omega(r)\to0$ as $r\downarrow0$ is not used.
 Recall that
$$ L_0 u = u_t + a^{jk} u_{x^j x^k}.
$$ Let $\bQ$ be the collection of all $Q_r(t,x)$, $(t,x) \in \bR^{d+1}$,
$r \in (0, \infty)$. For a function $g$ defined on $\bR^{d+1}$, we denote
its (parabolic) maximal and sharp function, respectively, by
$$ M g (t,x) = \sup_{(t,x) \in Q} \dashint_{Q} |g(s,y)| \, dy \, ds,
$$
$$ g^{\#} (t,x) = \sup_{(t,x) \in Q} \dashint_{Q} |g(s,y) - (g)_{Q}| \, dy
\, ds,
$$ where the supremums are taken over all $Q \in \bQ$ containing $(t,x)$.

\begin{theorem} \label{main_sharp} 
Let $\mu$, $\nu \in (1,\infty)$,
$1/\mu + 1/\nu = 1$, and $R \in (0, \infty)$. There exists a constant $N
= N(d, \delta, \mu)$ such that, for any $u \in C_0^{\infty}(Q_R)$, we have
$$
 (u_t)^{\#} + (u_{xx'})^{\#} \le N ( a_R^{\#} )^{\alpha/\nu} 
\left[ M ( |u_{xx}|^{2 \mu} ) \right]^{1/(2\mu)}
$$
\begin{equation}
                                          \label{3.13.3}
 + N \left[ M ( |L_0 u|^{2} ) \right]^{\alpha}
\left[ M ( |u_{xx}|^{2} ) \right]^{\beta} + N \left[ M ( |L_0 u|^{2} )
\right]^{1/2},
\end{equation}
 where $\alpha = 1/(d + 4)$ and $\beta = (d + 2)/(2d + 8)$.
\end{theorem}
\begin{proof} Let $\kappa \ge 4$, $r \in (0, \infty)$, and $(t_0, x_0) =
(t_0, x^1_0, x_0') \in \bR^{d+1}$. Also recall that
the sets $\Gamma_{r}(t,x')$ are introduced
in Section \ref{section 3.13.1}
and  set 
$$
\bar{a}^{jk}(x^1) = 
\dashint_{\Gamma_{\kappa r}(t_0, x'_0)} a^{jk}(s,x^1,y') \, dy' \, ds
\quad \text{if} \quad \kappa r < R,
$$
$$
\bar{a}^{jk}(x^1) = 
\dashint_{\Gamma_{R}} a^{jk}(s,x^1,y') \, dy' \, ds
\quad \text{if} \quad \kappa r \ge R.
$$ For $\rho > 0$, we denote 
$$
\cA_{\rho} = (|L_0 u|^2)_{Q_{\rho}(t_0, x_0)}, \quad 
\cB_{\rho} = (|u_{xx}|^2)_{Q_{\rho}(t_0, x_0)}, 
\quad
\cC_{\rho} = (|u_{xx}|^{2 \mu})_{Q_{\rho}(t_0, x_0)}^{1/\mu}.
$$

Set $\bar{L}_0 u = u_t + \bar{a}^{jk} u_{x^j x^k}$. Also set $w$ to be
either $u_t$ or $u_{xx'}$. Then by Lemma \ref{u_sharp}, we have
\begin{equation}\label{sharp_u}
\left( |w - (w)_{Q_r(t_0, x_0)}|^2 \right)_{Q_r(t_0, x_0)}
\le N \kappa^{d+2} \left( |\bar{L}_0 u|^2 \right)_{Q_{\kappa r}(t_0, x_0)}
+ N \kappa^{-2} \cB_{\kappa r}.
\end{equation} Note that
\begin{equation} \label{esti_chi}
\int_{Q_{\kappa r}(t_0, x_0)} |\bar{L}_0 u|^2 \, dx \, dt
\le 2 \int_{Q_{\kappa r}(t_0, x_0)} |L_0 u|^2 \, dx \, dt + I,
\end{equation} where $I$ is a constant times
$$
\int_{Q_{\kappa r}(t_0,x_0)} | ( \bar{L}_0 - L_0 ) u |^2 \, dx \, dt =
\int_{Q_{\kappa r}(t_0,x_0) \cap Q_R} \dots
\le J_1^{1/\nu} J_2^{1/\mu},
$$
$$ J_1 = \int_{Q_{\kappa r}(t_0, x_0) \cap Q_R} | \bar{a} - a |^{2 \nu}
\, dx \, dt,
$$
$$ J_2 = \int_{Q_{\kappa r}(t_0, x_0) \cap Q_R}  |u_{x x}|^{2 \mu} \, dx
\, dt \le N (\kappa r)^{d+2} (\cC_{\kappa r})^{\mu}.
$$
Observe that  if $\kappa r < R$,
$$ J_1
\le N \int_{x_0^1 - \kappa r}^{x_0^1 + \kappa r}\int_{\Gamma_{\kappa
r}(t_0, x'_0)} | \bar{a} - a | \, dx' \, dt \, dx^1
$$
$$
\le N (\kappa r)^{d+2} a^{\#}_{\kappa r}
\le N (\kappa r)^{d+2} a^{\#}_{R}.
$$ In case $\kappa r \ge R$,
$$ J_{1}
\le N \int_{-R}^{R}\int_{\Gamma_{R}} | \bar{a} - a | \, dx' \, dt \, dx^1
\le N R^{d+2} a^{\#}_{R}
\le N (\kappa r)^{d+2} a^{\#}_{R}.
$$ 
>From  \eqref{esti_chi} and the above estimates, we have
$$
\left( |\bar{L}_0 u|^2 \right)_{Q_{\kappa r}(t_0, x_0)}
\le N ( a^{\#}_R )^{1/\nu} \cC_{\kappa r} + N \cA_{\kappa r}.
$$ This, together with \eqref{sharp_u}, gives us
\begin{multline}\label{Q_r_ineq}
\left(|w - (w)_{Q_r(t_0, x_0)}|^2 \right)_{Q_r(t_0, x_0)} 
\le N \kappa^{d+2} (a^{\#}_R )^{1/\nu} \cC_{\kappa r} \\ + N \kappa^{d+2}
\cA_{\kappa r} + N \kappa^{-2} \cB_{\kappa r}.
\end{multline}

Now observe that $\cB_{\kappa r} \le M(|u_{xx}|^2)(t,x)$ for any $(t,x)
\in Q_r(t_0,x_0)$. Similar inequalities hold true for $\cA_{\kappa r}$
and $\cC_{\kappa r}$.   From this fact and  
\eqref{Q_r_ineq} it follows that, for any $(t,x) \in \bR^{d+1}$ and $Q
\in \bQ$ such that $(t,x) \in Q$,
$$
\left(|w - (w)_{Q}|^2 \right)_{Q} 
\le N \left( \kappa^{d+2} (a^{\#}_R )^{1/\nu} \cC (t,x) + \kappa^{d+2}
\cA (t,x) + \kappa^{-2} \cB (t,x) \right),
$$ where
$$
\cA = M(|L_0 u|^2), 
\quad 
\cB = M(|u_{xx}|^2), 
\quad
\cC = \left( M(|u_{xx}|^{2 \mu}) \right)^{1/\mu}.
$$

 Note that the above inequality is proved for $\kappa \ge 4$. In case $0
< \kappa < 4$, we have
$$
\dashint_{Q} |w - (w)_{Q}|^2 \, dx \, dt 
\le \left(|w|^2 \right)_{Q}
\le N \left( |L_0 u|^2 \right)_{Q} + N \left( |u_{xx}|^2 \right)_{Q} 
$$
$$
\le N \left( \cA(t,x) + 16 \kappa^{-2} \cB(t,x) \right)
$$ for $(t,x) \in Q \in \bQ$. Therefore, we finally have
$$
\left(|w - (w)_{Q}|^2 \right)_{Q} 
\le N \kappa^{d+2} (a^{\#}_R)^{1/\nu} \cC (t,x) 
$$
$$ + N (\kappa^{d+2}+1) \cA (t,x) + N \kappa^{-2} \cB (t,x)
$$ for all $\kappa > 0$, $(t,x) \in \bR^{d+1}$, and $Q \in \bQ$
satisfying $(t,x) \in Q$.

 Take the supremum of the left side of the above
inequality over all $Q \in \bQ$ containing $(t,x)$, and then minimize the
right-hand side with respect to 
$\kappa > 0$. 
Also observe that
$$
\big(\dashint_{Q} |w - (w)_{Q}| \, dx \, dt\big)^{2}\leq
\dashint_{Q} |w - (w)_{Q}|^2 \, dx \, dt
$$

Then we obtain
$$
[w^{\#}]^{2}(t,x)\leq N\cA(t,x)+
\big[(a^{\#}_R)^{1/\nu} \cC +\cA\big]^{2/(d+4)}
\cB^{(d+2)/(d+4)}(t,x).
$$
Here $\cB\leq\cC$ and this leads to \eqref{3.13.3}.\end{proof}

\begin{corollary} For $p>2$, there exist constants $R = R(d,\delta,
p,\omega)$ and $N = N(d,\delta,p)$ such that, for any $u \in
C_0^{\infty}(Q_R)$, we have
$$
\| u_{xx} \|_{L_p(\bR^{d+1})} \le N \|L_0 u\|_{L_p(\bR^{d+1})}.
$$
\end{corollary}

\begin{proof} Set $L_p = L_p(\bR^{d+1})$. Choose a number $\mu$ such that
$p > 2 \mu > 1$. Then we use \eqref{3.13.3}
together with the Fefferman-Stein theorem on sharp functions,
H\"{o}lder's inequality,  and the Hardy-Littlewood maximal function
theorem to obtain
\begin{equation}\label{L_p_esti}
\| u_t \|_{L_p} + \| u_{xx'} \|_{L_p}
\le N (a_R^{\#})^{\alpha/\nu} \| u_{xx} \|_{L_p} + N \|L_0
u\|_{L_p}^{2 \alpha} \|u_{xx}\|_{L_p}^{2 \beta} + N \|L_0 u\|_{L_p},
\end{equation}
 where, as noted in Theorem \ref{main_sharp},
$1/\mu + 1/\nu =1$ and $2 \alpha + 2 \beta = 1$. Now we notice that
$$ u_{x^1x^1} = \frac{1}{a^{11}} \left( L_0 u - u_t - \sum_{j \ne 1, k
\ne 1} a^{jk} u_{x^j x^k} \right).
$$ 
>From this and \eqref{L_p_esti}, we have
$$
\|u_{xx}\|_{L_p}
\le N (a_R^{\#})^{\alpha/\nu} \| u_{xx} \|_{L_p} + N \|L_0
u\|_{L_p}^{2 \alpha} \|u_{xx}\|_{L_p}^{2 \beta} + N \| L_0 u \|_{L_p}. 
$$ 
Choose an appropriate $R$ such that 
$$  N (a_R^{\#})^{\alpha/\nu} \le 1/2.
$$ Then 
$$
\frac{1}{2}\| u_{xx} \|_{L_p} 
\le  N \|L_0 u\|_{L_p}^{2 \alpha} \|u_{xx}\|_{L_p}^{2 \beta} + N \| L_0 u
\|_{L_p}. 
$$ This finishes the proof. ~\end{proof}

\begin{lemma}\label{ind_T} Let $T \in (0, \infty]$.  Then there exists
$\lambda_0 = \lambda_0(d,\delta,K,p,\omega) \ge 0$ such that, for all
$\lambda \ge \lambda_0$ and $u \in \WO{1,2}{p}(\Omega_T)$,
$$
\lambda \| u \|_{L_p(\Omega_T)} + \| u_{xx} \|_{L_p(\Omega_T)} + \| u_t
\|_{L_p(\Omega_T)}
\le N \| L u - \lambda u \|_{L_p(\Omega_T)},
$$ where $N = N(d,\delta,K,p,\omega)$ (independent of $T$).
\end{lemma}

\begin{proof} We have an $L_p$-estimate for functions with small compact
support. Thus the rest of the proof can be done by following the argument
in \cite{Kr05}. ~\end{proof}

Now Theorem ~\ref{para_main_01} follows from  the above lemma
and the argument in \cite{Kr05}. This ends the proof of Theorem
~\ref{para_main_01}.

\end{document}